\documentstyle{amsart}
\newtheorem{theorem}{Theorem}[section]

\theoremstyle{remark}
\newtheorem{remark}[theorem]{Remark}

\newcommand{\FF}{{\Bbb{F}}}
\newcommand{\Q}[2]{Q_{{#1},{#2}}}
\newcommand{\EQ}[3]{Q^{#1}_{{#2},{#3}}}
\newcommand{\sq}[1]{{\rm Sq}^{#1}}
\newcommand{\tk}{\,\,the $\chi$-{\em trick}\,\,}

\begin{document}
\title[Dickson Invariants in the image of the Steenrod Square
]{Dickson Invariants in the image of the Steenrod Square
}
\author[Kai Xu ]{Kai Xu}
\email[Kai Xu]{matxukai@@leonis.nus.edu.sg}
\address{Department of Mathematics, National University of Singapore, 2 Science Drive 2,
117543 Singapore}
\date{}
\begin{abstract}
Let $D_n$ be the Dickson invariant ring of $\FF_2[X_1,\ldots,X_n]$ acted by the general linear
group $\text{GL}(n,\FF_2)$. 
In this paper, we provide an elementary proof of the  conjecture by 
\cite{hung2}:
each element in $D_n$ is in the image of  the Steenrod square in
$\FF_2[X_1,\ldots,X_n]$, where $n>3$. 
\end{abstract}
\maketitle
\section{Introduction}

A polynomial in $\FF_2[X_1,X_2,\ldots, X_n]$ is {\em hit} if it is in the image of
the summation of the Steenrod square: 
$\sum_{i\geq 1}\sq{i}$. 
Let $D_n$ be the Dickson invariant algebra of $n$-variables.
In this paper, we will prove the following,

\begin{theorem}\label{main}
When $n>3$, each polynomial in the Dickson invariant ring $D_n$ is hit.
\end{theorem}

In \cite{hung2}, Hung studies the Dickson invariants in the image of the Steenrod
square. Since it is trivial that $D_1$ and $D_2$ are not hit,
the problem starts interesting from $n=3$. In the same paper, Hung shows that
each element in $D_3$ is hit and conjectured  that it is true for $D_{n>3}$.
So our result provides a positive answer to the conjecture, which supports to the 
positive answer of the conjecture on the spherical classes: there are no spherical classes
in $Q_0S^0$, except the Hopf invariant one and Kervaire invariant one elements.
We refer to \cite{hung2} and 
an excellent expository paper~\cite{wood2}, p501 for more background
regarding to this conjecture. 

\begin{remark}
Recently, K. F. Tan and the author \cite{TX} has obtained an  elementary proof of 
the case $n=3$. 
\end{remark}

\section{Proof of Theorem \ref{main}}

We first recall some basic properties regarding the Dickson algebra.

Write $V_n$ for the product
$$
\prod_{\alpha_i\in\{0,1\},i=1,..,n-1}(\alpha_1x_1+\cdots+\alpha_{n-1}x_{n-1}+x_n).
$$
Then we have the following theorem.
\begin{theorem}[Hung~\cite{hung1}]\label{st}
$$\sq{i}V_n=\begin{cases}
V_n \text{ if $i=0$}\\
V_n\EQ{}{n-1}{s} \text{ if $i=2^{n-1}-2^s$, $0\leq s\leq n-1$}\\
V^2_n \text{ if $i=2^{n-1}$}\\
0  \text{ otherwise.}
\end{cases}
$$
$$\sq{i}\EQ{}{n}{s}=\begin{cases} \EQ{}{n}{r} \text{ if $i=2^s-2^r$, $r\leq s$}\\
\EQ{}{n}{r}\EQ{}{n}{t} \text{ if $i=2^n-2^t+2^s-2^r$, $r\leq s <t$}\\
\EQ{2}{n}{s} \text{ if $i=2^n-2^s$ }\\
0 \text{ otherwise.}
\end{cases}
$$
\end{theorem}
In the following, we will frequently use the above results without mentioning 
each time.

We use the induction on $n$ to prove Theorem \ref{main}. 
Suppose that the statement is true for $n$. Then we will prove that each polynomial
in $D_{n+1}$ is hit.

Recall that 
$$\EQ{}{n+1}{k}=\EQ{2}{n}{k-1}+V_{n+1}\EQ{}{n}{k}\quad \text{ for $1\leq k\leq n.$}$$
So any monomial in $\FF_2[\Q{n+1}{0}, \Q{n+1}{1},... ,\Q{n+1}{n}]$ 
can be written as the summation of  the following form:
$$
A:=V^a_{n+1}\EQ{n_0}{n}{0}\EQ{n_1}{n}{1}\EQ{n_2}{n}{2}\cdots\EQ{n_{n-1}}{n}{n-1}.
$$
Hence by the hypothesis of the induction,
it is sufficient to show that  $A$ is hit for any $a>0$.
Notice that 
\begin{equation}\label{V}
V_{n+1}=\sum^{n}_{s=1}\sq{1}(\EQ{}{n}{s}X^{2^s-1}_{n+1}).
\end{equation}

When $n_1$ is even, we have the hit polynomial
$$A=\sq{1}[\left(\sum^{n}_{s=1}\EQ{}{n}{s}x^{2^s-1}_{n+1}\right)
V^{a-1}_{n+1}\EQ{n_0}{n}{0}\EQ{n_1}{n}{1}\EQ{n_2}{n}{2}\cdots\EQ{n_{n-1}}{n}{n-1}].
$$

If $n_1$ is odd and $n_2$ is  even, then $A$ can be 
written as the hit polynomial:
\begin{multline*}
\sq{2}[V^{a}_{n+1}\EQ{n_0}{n}{0}\EQ{n_1-1}{n}{1}\EQ{n_2+1}{n}{2}\cdots\EQ{n_{n-1}}{n}{n-1}]\\
+\sq{1}[\left(\sum^{n}_{s=1}\EQ{}{n}{s}x^{2^s-1}_{n+1}\right)V^{a-1}_{n+1}\EQ{n_0}{n}{0}(\sq{1}
\EQ{\frac{n_1-1}{2}}{n}{1})^2\EQ{n_2+1}{n}{2}\cdots\EQ{n_{n-1}}{n}{n-1}]
\end{multline*}

In the following, we will always assume that $n_1$ and $n_2$ are both odd.

When $n=3$, $n_0$ is even and $a$ is odd, we have 
$$
\begin{array}{lll}
A&=&(V^{a-1}_{4}\sq{4}V_{4})\EQ{n_0}{3}{0}\EQ{n_1}{3}{1}\EQ{n_2-1}{3}{2}\\
&\equiv&
V_{4}\chi(\sq{4})[V^{a-1}_{4}\EQ{n_0}{3}{0}\EQ{n_1}{3}{1}\EQ{n_2-1}{3}{2}]
\text{ (modulo the hits)}\\
&\equiv& V^{a}_4\EQ{}{3}{1}\left(\sq{2}[\EQ{\frac{n_0}{2}}{3}{0}\EQ{\frac{n_1-1}{2}}{3}{1}
\EQ{\frac{n_2-1}{2}}{3}{2}]\right)^2\text{ (modulo the hits)}
.
\end{array}
$$
Using the previous observation,
the last polynomial is hit, since the order of $\EQ{}{3}{2}$ is even.

When $n=3$, $n_0$ is even and $a$ is even,
notice that 
$$\EQ{n_0}{3}{0}\EQ{n_1}{3}{1}\EQ{n_2}{3}{2}=
\EQ{n_0}{3}{0}\EQ{n_1-1}{3}{1}\EQ{n_2-1}{3}{2}\sq{4}\EQ{}{3}{1}.$$
Then using \tk and doing some basic computation, we can see that the monomial
$\EQ{n_0}{3}{0}\EQ{n_1}{3}{1}\EQ{n_2}{3}{2}$ is 
in the image of $\sum^{4}_{i=1}\sq{i}$.
In fact, 
$$
\begin{array}{lll}
\EQ{}{3}{1}\chi(\sq{4})[\EQ{n_0}{3}{0}\EQ{n_1-1}{3}{1}\EQ{n_2-1}{3}{2}]&&\\
=
[\sq{2}\EQ{}{3}{2}][\sq{2}
(\EQ{\frac{n_0}{2}}{3}{0}\EQ{\frac{n_1-1}{2}}{3}{1}\EQ{\frac{n_2-1}{2}}{3}{2})]^2&&\\
\equiv  \EQ{}{3}{2} \chi(\sq{2})[\EQ{\frac{n_0}{2}}{3}{0}\EQ{\frac{n_1-1}{2}}{3}{1}\EQ{\frac{n_2-1}{2}}{3}{2}]^2
\text{ (modulo the hits)} &&\\
=(\EQ{2}{2}{1}+V_3)[\sq{1}\sq{2}(\EQ{\frac{n_0}{2}}{3}{0}\EQ{\frac{n_1-1}{2}}{3}{1}\EQ{\frac{n_2-1}{2}}{3}{2})]&&\\
= \sq{2}\left(\EQ{}{2}{1}[\sq{1}\sq{2}(
\EQ{\frac{n_0}{2}}{3}{0}\EQ{\frac{n_1-1}{2}}{3}{1}\EQ{\frac{n_2-1}{2}}{3}{2})^2]\right) &&\\
+
\sq{1}\left\{(\EQ{}{2}{1}X_3+X^3_3)[\sq{1}\sq{2}(\EQ{\frac{n_0}{2}}{3}{0}\EQ{\frac{n_1-1}{2}}{3}{1}\EQ{\frac{n_2-1}{2}}{3}{2})]^2\right\},&&
\end{array}
$$
where we have used (\ref{V}) in the last equality.
On the other hand, 
$\sq{i}V^{a}_{4}=0$ for $i=1,2,3$ and $4$. Therefore using  \tk,
we know that $A$ is hit.

When $n=3$, $n_0$ is odd and $a$ is odd, the polynomial $A$ equals 
$$
V^{a-1}_4(\sq{7}V_{4})\EQ{n_0-1}{3}{0}\EQ{n_1}{3}{1}\EQ{n_2}{3}{2}
.
$$
From  the discussion above, we know that
$$\EQ{n_0-1}{3}{0}\EQ{n_1}{3}{1}\EQ{n_2}{3}{2}$$ is 
in the image of $\sum^{4}_{i=1}\sq{i}$. On the other hand,
$\chi(\sq{i})(V^{a-1}_4\sq{7}V_{4})=0$  for $i=1,2,3$ and $4$. So using  \tk,
we conclude that $A$ is hit.

When $n=3$, $n_0$ is odd and $a$ is even, let $\nu$ be the integer such that
$a=2^{\nu}b$ where $b$ is odd. Then 
$$
V^{a}_4=\sq{4a}\sq{2a}\cdots\sq{8b}V^{b}_4.
$$
Hence 
$$
\begin{array}{lll}
A&=&(\sq{4a}\sq{2a}\cdots\sq{8b}V^{b}_4)(\EQ{n_0}{3}{0}\EQ{n_1}{3}{1}\EQ{n_2}{3}{2})\\
&\equiv& V^{b}_4\chi(\sq{8b})\cdots\chi(\sq{2a})\chi(\sq{4a})
(\EQ{n_0}{3}{0}\EQ{n_1}{3}{1}\EQ{n_2}{3}{2})\,\,\, \text{ (modulo the hits)}.
\end{array}
$$
After expanding the last polynomial using Theorem \ref{st}, it is easy to see that 
each resulting term belongs to one of the previous cases. Therefore $A$ is hit.

When $n\geq 4$, the polynomial $A$ takes the following  form,
\begin{equation}\label{ol}
V^{a}_{n+1}(\sq{2^n-4}\EQ{}{n}{1})\EQ{n_0}{n}{0}\EQ{n_1-1}{n}{1}\EQ{n_2-1}{n}{2}
\cdots\EQ{n_{n-1}}{n}{n-1}.
\end{equation}
Using a result of Don Davis, Theorem 2. of \cite{DD} and \tk, we know that 
it is sufficient to show the polynomial:
$$
\EQ{}{n}{1}\sq{2^{n-1}}\cdots\sq{8}\chi(\sq{4})\left\{V^{a}_{n+1}
\EQ{n_0}{n}{0}\EQ{n_1-1}{n}{1}\EQ{n_2-1}{n}{2}\cdots\EQ{n_{n-1}}{n}{n-1}\right\}
$$
is hit.

After expansion using the Steenrod operation, the above polynomial can be written as the summation of the form:
$$
V^{a}_{n+1}\EQ{k_0}{n}{0}\EQ{k_1}{n}{1}\EQ{k_2}{n}{2}\cdots\EQ{k_{n-1}}{n}{n-1}.
$$
Using the previous discussion, we can conclude that all these polynomials are hit,
except for those when $k_1$ and $k_2$ are both odd. But in this case, we can replace
$n_i$ by $k_i$ for all $i$ in (\ref{ol})
 and carry out the above process again.
After using this process sufficiently many times
with modulo the hits, we can conclude that the new  $k_0$, $k_1$  and $k_2$ 
are independent of the process.
To keep $k_0$, $k_1$ and $k_2$ unchanged with the process,
we must require that 
$$
\sq{2^{n-1}}\cdots\sq{8}\chi(\sq{4})\left\{V^{a}_{n+1}
\EQ{k_3}{n}{3}\cdots\EQ{k_{n-1}}{n}{n-1}\right\}  \text{ (modulo the hits)}
$$
contributes $\EQ{}{n}{2}$ after each process is done, since for $j\leq 2^{n-1}$ and $t<n$ ,
$\sq{j}\EQ{}{n}{0}=\EQ{}{n}{0}\EQ{}{n}{t}$ only if $j=t=n-1(>2)$.
Finally because all $k_i$ ($0\leq i<n$)
are finite, we conclude that  $A$ is hit after carrying on the process further for
enough many times.

\end{document}